\documentclass{article}

\usepackage{amssymb}
\usepackage{amsmath}
\usepackage{amsfonts}
\usepackage{amscd}

\newtheorem{theo}{Th\'{e}or\`{e}me}
\newtheorem{defi}{D\'{e}finition}
\newtheorem{propo}{Proposition}

\newtheorem{coro}{Corollaire}

\usepackage[dvips]{graphicx}
\newcommand{\C}{\mathbb C}
\newcommand{\g}{\frak{g}}
\newcommand{\w}{\frak{t}}
\newcommand{\h}{\frak{h}}

\newcommand{\n}{\frak{n}}

\newcommand\dd{\noindent{\it D\'emonstration. }}

\newcommand{\al}{alg\`ebre de Lie }
\newcommand{\als}{alg\`ebres de Lie }

\begin{document}

\title{Sur les alg\`ebres de Lie quasi-filiformes compl\'etables}
\author{Luc\'{\i}a Garc\'{\i}a Vergnolle
\footnote{Dpto. Geometr\'{\i}a y Topolog\'{\i}a, Facultad de
Ciencias Matem\'aticas U.C.M. Plaza de Ciencias 3, 28040 Madrid,
Espagne lucigarcia@mat.ucm.es} }

\date{}

\maketitle

\begin{abstract}
Le but de ce travail est  de d\'eterminer les alg\`ebres
quasi-filiformes compl\'etables. Nous prouvons de plus que, pour
tout entier positif $m$, il existe une \al compl\`ete dont la
dimension du deuxi\`eme groupe de cohomologie est sup\'erieure ou
\'egale \`a $m$.
\end{abstract}

{\small Mots clefs : compl\`ete, compl\'etable, quasi-filforme, \al}

\section{Introduction}

La compl\'etude d'une alg\`ebre de Lie, \'etant une propri\'et\'e
d\'eduite de la structure des d\'erivations, constitue un
invariant d'int\'er\^{e}t pour l'\'etude du comportement d'une
classe d'isomorphisme par rapport aux d\'eformations et
contractions. Bien que la notion d'\al compl\`ete e\^ut \'et\'e
introduite en 1951 dans le contexte de la th\'eorie d'alg\`ebres
sous-invariantes de Schenckman ~\cite{Sch}, Jacobson fut le
premier \`a donner une d\'efinition formelle dans les ann\'ees 60,
en utilisant des outils
cohomologiques.\\
Par la suite, de nombreux auteurs se sont int\'er\'ess\'es \`a
l'\'etude des alg\`ebres de Lie compl\`etes. Favre a \'etudi\'e
les \als compl\`etes par rapport \`a leur  nilradical ~\cite{MFa},
tandis que Carles a analys\'e la suite croissante des \als de
d\'erivations ~\cite{Ca}. Plus r\'ecemment, Zhu et Meng ont
\'etudi\'e la compl\'etude des \als r\'esolubles de rang maximal
~\cite{Zhu} et non-maximal ~\cite{Meng}.
\begin{defi}
Une \al $\g$ est dite compl\`ete si
\begin{enumerate}
\item le centre de $\g$ est nul, $Z(\frak{g})=\{0\}$
\item toutes ses d\'erivations sont int\'erieures, c'est-\`a-dire, $Der(\frak{g})=ad(\frak{g})$.
\end{enumerate}
\end{defi}
Soit $\g$ une \al et $H^{n}(\g,\g)$ son n-i\`eme groupe de cohomologie. Rappelons que $H^{0}(\g,\g)=Z(\g)$. De plus, l'ensemble des cocycles $Z^{1}(\frak{g},\frak{g})$ correspond \`a l'ensemble des d\'erivations de $\frak{g}$ et l'ensemble des cobords $B^{1}(\frak{g},\frak{g})$ correspond \`a l'ensemble des d\'erivations int\'erieures. On en d\'eduit que l'\al $\g$ est compl\`ete si et seulement si $H^{0}(\g,\g)=H^{1}(\g,\g)=\{0\}$ .\\
L'\'etude du deuxi\`eme groupe de cohomologie est li\'ee \`a celle des d\'eformations. Si l'\al $\g$ est donn\'ee par les crochets $[,]$, une d\'eformation formelle de $\g$ est d\'efinie par la s\'erie formelle:
\begin{equation*}
\phi_{t}(X,Y)=[X,Y]_{t}=F_{0}(X,Y)+F_{1}(X,Y)t+\dots, \quad \forall X,Y \in \g
\end{equation*}
avec $F_{0}(X,Y)=[X,Y], \quad \forall X,Y \in \g$.\\
En imposant la condition de Jacobi sur $[,]_{t}$, on d\'eduit que $F_{1}\in Z^{2}(\frak{g},\frak{g})$. Ainsi, \`a chaque d\'eformation on peut faire correspondre un $2$-cocycle. Les d\'eformations infinit\'esimales sont celles qui v\'erifient la condition de Jacobi jusqu'au premier ordre et s'identifient aux \'el\'ements de $Z^{2}(\frak{g},\frak{g})$.\\
Par ailleurs, on dit que deux d\'eformations $\phi_{t}$ et $\phi^{\prime}_{t}$ de $\g$ sont \'equivalentes s'il existe une s\'erie formelle $g(t)=\sum_{p=0}^{\infty}G_{p}t^{p}$ avec $G_{p}\in GL(n,\mathbb{C})$  telle que
\begin{equation*}
\phi_{t}(X,Y)= (g(t)\phi^{\prime}_{t})(g(t)^{-1}X,g(t)^{-1}Y),\quad \forall X,Y \in \g.
\end{equation*}
Si deux d\'eformations sont \'equivalentes alors les $2$-cocycles correspondants sont \'egaux modulo $B^{2}(\frak{g},\frak{g})$, la r\'eciproque \'etant vraie aussi pour les d\'eformations lin\'eaires ~\cite{Vi}.
\par
Soit $\n$ une \al nilpotente de dimension finie sur $\C$ et $Der(\n)$ son alg\`ebre de d\'erivations. Un tore $\w$ sur $\n$ est une sous-alg\`ebre commutative de $Der(\n)$ form\'ee par des endomorphismes semi-simples. Il est clair que $\n$ se d\'ecompose de la fa\c{c}on suivante
\begin{equation*}
\n=\sum_{\alpha \in \w^{\ast}} \n_{\alpha}
\end{equation*}
o\`u $\w^{\ast}$ repr\'esente l'espace dual de $\w$ et $\n_{\alpha}=\{ X\in \n \;|\; [f,X]=\alpha(f)X, \quad \forall f\in \w \}$. Lorsque $\w$ est un tore maximal par rapport \`a l'inclusion, l'ensemble $\Delta=\{\alpha\in\w^{\ast} \,|\, \n_{\alpha}>0\}$ est un syst\`eme de poids de $\n$ ~\cite{Fa}. Si $\w$ et $\w^{\prime}$ sont deux tores maximaux, ils sont conjugu\'es par automorphismes, c'est-\`a-dire, il existe un automorphisme $\theta\in Aut(\n)$ tel que $\theta \w \theta^{-1}= \w^{\prime}$. Le rang d'un tore maximal est donc un invariant de $\n$ que l'on appelle le rang de $\n$. Par ailleurs, le type de $\n$, not\'e par $ty(\n)$, est la dimension du quotient $H^{1}(\g,\C)=\n/[\n,\n]$. On montre que le rang de $\n$ est toujours major\'e par son type ~\cite{Fa}. L'\al $\n$ est de rang maximal si le rang et le type sont \'egaux.\\
Soit $\w$ un tore maximal , on d\'efinie $\g=\w\overrightarrow{\oplus}\n$ comme:
$$
[f_{1}+x_{1},f_{2}+x_{2}]=f_{1}(x_{2})-f_{2}(x_{1})+[x_{1},x_{2}] \quad \forall f_{1},f_{2} \in \w,\, \forall x_{1},x_{2} \in \n
$$
L' \al $\g$ est alors r\'esoluble et $\w$ est une sous-alg\`ebre de Cartan de $\g$. On dira que le rang de $\g$ est celui de $\n$ et que $\g$ est de rang maximal si $\n$ l'est aussi.
\begin{theo}
~\cite{Zhu}
Si $\g$ est une \al r\'esoluble compl\`ete, elle se d\'ecompose de la fa\c{c}on $\h\overrightarrow{\oplus}\n$, $\n$ \'etant le nilradical et $\h$ une sous-alg\`ebre isomorphe \`a un tore maximal de $\n$. De plus, $\h$ est une sous-alg\`ebre de Cartan de $\g$.
\end{theo}
\begin{defi}
Une alg\`ebre de Lie nilpotente $\n$ est compl\'etable lorsque la somme semi-directe $\h\overrightarrow{\oplus}\n$, $\h$ \'etant un tore maximal de $\n$, est compl\`ete.
\end{defi}
Les alg\`ebres nilpotentes les plus \'etudi\'ees sont les filiformes.

\section{Alg\`ebres de Lie quasi-filiformes compl\'etables}

Soit $\frak{g}$ une alg\`{e}bre de Lie nilpotente de dimension $n$ et nilindice $m$. Elle est naturellement filtr\'{e}e par la suite centrale descendante:\\
\begin{align*}
\frak{g}_{1}  &  =\frak{g}\supseteq \frak{g}_{2}=\left[
\frak{g},\frak{g}\right]\supseteq \frak{g}_{3}=\left[
\frak{g}_{2},\frak{g}\right]  \supseteq...\supseteq
\frak{g}_{k+1}=\left[  \frak{g}_{k},\frak{g}\right]  \supseteq...
\supseteq \frak{g}_{m+1}=\left\{ 0 \right\}
\end{align*}
On peut alors associer une alg\`{e}bre de Lie gradu\'{e}e, not\'{e}e par $gr(\frak{g})$, et d\'{e}finie par:
$$
{\rm gr}\frak{g} \; = \; \sum_{i=1}^{m} \, \frac{\frak{g}_{i}}{\frak{g}_{i+1}} \; = \; \sum_{i=1}^{m} W_{i},
$$
dont le crochet est donn\'{e} par:
$$[X+\frak{g}_{i+1},Y+\frak{g}_{j+1}]=[X,Y]+\frak{g}_{i+j+1},\quad \forall X \in \frak{g}_{i},\quad \forall Y \in \frak{g}_{j}.$$
Si on consid\`ere la suite $\{p_{1},\dots,p_{m}\}$ o\`u $p_{i}=\dim \frac{\frak{g}_{i}}{\frak{g}_{i+1}}=\dim W_{i}$, cette suite est la m\^eme pour $gr(\frak{g})$ que pour $\g$, et en particulier, $p_{1}=ty(gr(\g))=ty(\g)$.
\begin{defi}
Une alg\`{e}bre $\frak{g}$ est gradu\'{e}e naturellement quand elle est isomorphe \`{a} ${\rm gr}\,\frak{g}$.
\end{defi}
Les alg\`{e}bres filiformes sont celles dont le nilindice est maximal. Il est clair que l'alg\`{e}bre gradu\'{e}e d'une alg\`{e}bre filiforme est aussi filiforme. La classification des alg\`{e}bres de Lie filiformes gradu\'{e}es naturellement est due \`{a} Vergne ~\cite{Ver}.
Dans ~\cite{Go}, on \'{e}tudie une classe plus large d'alg\`{e}bres filiformes gradu\'{e}es, la graduation \'{e}tant d\'{e}finie par les racines d'un tore externe non nul de dimension quelconque.
A partir de cette classification et en utilisant les r\'esultats des articles ~\cite{Zhu} et ~\cite{Meng}, on prouve dans ~\cite{AnCa} le th\'eor\`eme suivant.
\begin{theo}
Toute \al filiforme de rang non nul est compl\'etable.
\end{theo}
On peut alors tenter de g\'en\'eraliser ce r\'esultat aux alg\`ebres dont le nilindice $m$ est \'egal \`a $n-2$ appel\'ees quasi-filiformes. Pour cela, faisons r\'ef\'erence \`a la classification des alg\`ebres quasi-filiformes naturellement gradu\'ees.
\begin{theo}
\cite{Gom} Soit $\frak{g}$ une alg\`{e}bre de Lie quasi-filiforme gradu\'{e}e naturellement de dimension $n$. Il existe alors une base $\{X_{0},X_{1},X_{2},\dots,X_{n-1}\}$ de $\frak{g}$ dans laquelle $\frak{g}$ est une des alg\`{e}bres d\'{e}crites ci-dessous.
\begin{enumerate}
\item Pour $\{p_{1}=3,p_{2}=1,p_{3}=1,\dots,p_{n-2}=1\}$
\begin{enumerate}
\item $L_{n-1}\oplus \mathbb{C} \quad (n\ge 4)$
$$
\lbrack X_{0},X_{i}\rbrack=X_{i+1}, \quad 1 \le i \le n-3.
$$
\item $Q_{n-1}\oplus \mathbb{C} \quad (n\ge 7,\: n\;{\rm impair})$
$$
\begin{array}{ll}
\lbrack X_{0},X_{i}\rbrack=X_{i+1}, & 1 \le i \le n-4,\\
\lbrack X_{i},X_{n-i-2}\rbrack=(-1)^{i-1}X_{n-2}, & 1\leq i\leq \frac{n-3}{2}.\\
\end{array}
$$
\end{enumerate}
\item Pour $\{p_{1}=2,p_{2}=1,\dots,p_{r-1}=1,p_{r}=2,p_{r+1}=1,\dots,p_{n-2}=1\}$ o\`{u} $r\in \{2,\dots,n-2\}$
\begin{enumerate}
\item $\frak{L_{n,r}}; \quad n\ge5, \;r\: {\rm impair}, \;3\le r\le 2[\frac{n-1}{2} ]-1$
$$
\begin{array}{ll}
\lbrack X_{0},X_{i} \rbrack=X_{i+1}, &  i=1,\dots,n-3\\
\lbrack X_{i},X_{r-i} \rbrack=(-1)^{i-1}X_{n-1}, &  i=1,\dots,\frac{r-1}{2}\\
\end{array}
$$
\item $\frak{Q_{n,r}}; \quad n\ge7,\;n\: {\rm impair}, \;r\: {\rm impair}, \;3\le r\le n-4$
$$
\begin{array}{ll}
\lbrack X_{0},X_{i} \rbrack=X_{i+1}, &  i=1,\dots,n-4\\
\lbrack X_{i},X_{r-i} \rbrack=(-1)^{i-1}X_{n-1}, &  i=1,\dots,\frac{r-1}{2}\\
\lbrack X_{i},X_{n-2-i} \rbrack=(-1)^{i-1}X_{n-2}, &  i=1,\dots,\frac{n-3}{2}\\
\end{array}
$$
\item $\frak{T_{n,n-4}};\quad n\ge7,\:n\; {\rm impair}$
$$
\begin{array}{ll}
\lbrack X_{0},X_{i} \rbrack=X_{i+1}, &  i=1,\dots,n-5\\
\lbrack X_{0},X_{n-3} \rbrack=X_{n-2},\\
\lbrack X_{0},X_{n-1} \rbrack=X_{n-3},\\
\lbrack X_{i},X_{n-4-i} \rbrack=(-1)^{i-1}X_{n-1}, &  i=1,\dots,\frac{n-5}{2}\\
\lbrack X_{i},X_{n-3-i} \rbrack=(-1)^{i-1} \frac{n-3-2i}{2} X_{n-3}, &  i=1,\dots,\frac{n-5}{2}\\
\lbrack X_{i},X_{n-2-i} \rbrack=(-1)^{i} (i-1) \frac{n-3-i}{2} X_{n-2}, &  i=2,\dots,\frac{n-3}{2}\\
\end{array}
$$
\item $\frak{T_{n,n-3}};\quad n\ge6,\:n \;{\rm pair}$
$$
\begin{array}{ll}
\lbrack X_{0},X_{i} \rbrack=X_{i+1}, &  i=1,\dots,n-4\\
\lbrack X_{0},X_{n-1} \rbrack=X_{n-2},\\
\lbrack X_{i},X_{n-3-i} \rbrack=(-1)^{i-1}X_{n-1}, &  i=1,\dots,\frac{n-4}{2}\\
\lbrack X_{i},X_{n-2-i} \rbrack=(-1)^{i-1} \frac{n-2-2i}{2} X_{n-2}, &  i=1,\dots,\frac{n-4}{2}\\
\end{array}
$$
\item $\frak{E_{9,5}^{1}}$
$$
\begin{array}{lll}
\lbrack X_{0},X_{i}\rbrack =X_{i+1}, & i=1,2,3,4,5,6, &\lbrack X_{0},X_{8}\rbrack =X_{6},\\
\lbrack X_{1},X_{4}\rbrack =X_{8}, & \lbrack X_{1},X_{5}\rbrack =2X_{6}, & \lbrack X_{1},X_{6}\rbrack =3X_{7},\\
\lbrack X_{2},X_{3}\rbrack =-X_{8}, & \lbrack X_{2},X_{4}\rbrack =-X_{6}, & \lbrack X_{2},X_{5}\rbrack =-X_{7},\\
\lbrack X_{2},X_{8}\rbrack =-3X_{7}. & &\\
\end{array}
$$
\item $\frak{E_{9,5}^{2}}$
$$
\begin{array}{lll}
\lbrack X_{0},X_{i}\rbrack =X_{i+1}, &i=1,2,3,4,5,6, & \lbrack X_{0},X_{8}\rbrack =X_{6},\\
\lbrack X_{1},X_{4}\rbrack =X_{8}, & \lbrack X_{1},X_{5}\rbrack =2X_{6}, & \lbrack X_{1},X_{6}\rbrack =X_{7},\\
\lbrack X_{2},X_{3}\rbrack =-X_{8}, & \lbrack X_{2},X_{4}\rbrack =-X_{6}, & \lbrack X_{2},X_{5}\rbrack =X_{7},\\
\lbrack X_{2},X_{8}\rbrack =-X_{7}, & \lbrack X_{3},X_{4}\rbrack =-2X_{7}. &\\
\end{array}
$$
\item $\frak{E_{9,5}^{3}}$
$$
\begin{array}{lll}
\lbrack X_{0},X_{i}\rbrack =X_{i+1}, &i=1,2,3,4,5,6, & \lbrack X_{0},X_{8}\rbrack =X_{6},\\
\lbrack X_{1},X_{4}\rbrack =X_{8}, & \lbrack X_{1},X_{5}\rbrack =2X_{6}, & \lbrack X_{2},X_{3}\rbrack =-X_{8},\\
\lbrack X_{2},X_{4}\rbrack =-X_{6}, & \lbrack X_{2},X_{5}\rbrack =2X_{7}, & \lbrack X_{3},X_{4}\rbrack =-3X_{7}.\\
\end{array}\\
$$
\item $\frak{E_{7,3}}$
$$
\begin{array}{lll}
\lbrack X_{0},X_{i}\rbrack =X_{i+1}, &i=1,2,3,4, & \lbrack X_{0},X_{6}\rbrack =X_{4},\\
\lbrack X_{1},X_{2}\rbrack =X_{6}, & \lbrack X_{1},X_{3}\rbrack =X_{4}, & \lbrack X_{1},X_{4}\rbrack =X_{5},\\
\lbrack X_{2},X_{6}\rbrack =-X_{5}. & &\\
\end{array}\\
$$
\end{enumerate}
\end{enumerate}

\end{theo}
\bigskip
En utilisant ce r\'esultat, on peut d\'eterminer les alg\`ebres quasi-filiformes ayant un tore non nul ~\cite{Ga}.
\begin{theo}
\label{theo}
Soit $\frak{g}$ une alg\`{e}bre de Lie quasi-filiforme de dimension $n$ admettant une d\'{e}rivation diagonale $f$ non nulle. Il existe alors une base de $\frak{g}$, $\{Y_{0},Y_{1},\dots,Y_{n-1}\}$, form\'{e}e de vecteurs propres de $f$ dont les crochets v\'{e}rifient l'un des cas suivants:
\begin{enumerate}

\item Si ${\rm gr}\frak{g}\simeq \frak{L_{n-1}}\oplus \mathbb{C} \quad(n\ge 4)$ alors $ty(\g)=3$ et
\begin{enumerate}
\item $\frak{g}=L_{n-1}\oplus \mathbb{C}$
$$
\lbrack Y_{0},Y_{i}\rbrack=Y_{i+1}, \quad 1 \le i \le n-3,
$$
$$
f\sim diag(\lambda_{0},\lambda_{1},\lambda_{0}+\lambda_{1},2\lambda_{0}+\lambda_{1},\dots,(n-3)\lambda_{0}+\lambda_{1},\lambda_{n-1}),
$$
$$
rang(\g)=3.
$$

\item \label{1.(b)} $\frak{g}=A_{n-1}^{k}(\alpha_{1},\dots,\alpha_{t-1})\oplus \mathbb{C}, \quad t=[\frac{n-k}{2}], \quad 2 \le k \le n-4$
$$
\begin{array}{lll}
\lbrack Y_{0},Y_{i}\rbrack=Y_{i+1}, & 1 \le i \le n-3, &\\
\lbrack Y_{i},Y_{i+1}\rbrack=\alpha_{i} Y_{2i+k}, & 1 \le i \le t-1, &\\
\lbrack Y_{i},Y_{j}\rbrack=a_{i,j} Y_{i+j+k-1}, & 1 \le i <j, & i+j \le n-k-1,\\
\end{array}
$$
$$
f\sim diag(\lambda_{0},k\lambda_{0},(k+1)\lambda_{0},(k+2)\lambda_{0},\dots,(k+n-3)\lambda_{0},\lambda_{n-1}),
$$
$$
rang(\g)=2.
$$
\item \label{1.(c)} $\frak{g}=L_{n-1} \overrightarrow{\oplus}_{l} \mathbb{C} \quad (2 \le l \le n-3)$
$$
\begin{array}{lll}
\lbrack Y_{0},Y_{i}\rbrack & =Y_{i+1}, & 1 \le i \le n-3,\\
\lbrack Y_{i},Y_{n-1}\rbrack & =Y_{i+l}, & 1 \le i \le n-l-2,\\
\end{array}
$$
$$
f\sim diag(\lambda_{0},\lambda_{1},\lambda_{0}+\lambda_{1},2\lambda_{0}+\lambda_{1},\dots,(n-3)\lambda_{0}+\lambda_{1},l\lambda_{0}),
$$
$$
rang(\g)=2.
$$
\item \label{1.(d)} $\frak{g}=A_{n-1}^{k}(\alpha_{1},\dots,\alpha_{t-1})\overrightarrow{\oplus}_{l} \mathbb{C} \quad t=[\frac{n-k}{2}] \quad 2 \le k \le n-4 \quad 2 \le l \le n-3$
$$
\begin{array}{lll}
\lbrack Y_{0},Y_{i}\rbrack=Y_{i+1}, & 1 \le i \le n-3, &\\
\lbrack Y_{i},Y_{n-1}\rbrack=Y_{i+l}, & 1 \le i \le n-l-2, &\\
\lbrack Y_{i},Y_{i+1}\rbrack=\alpha_{i} Y_{2i+k}, & 1 \le i \le t-1, &\\
\lbrack Y_{i},Y_{j}\rbrack=a_{i,j} Y_{i+j+k-1}, & 1 \le i <j, & i+j \le n-k-1,\\
\end{array}
$$
$$
f\sim diag(\lambda_{0},k\lambda_{0},(k+1)\lambda_{0},(k+2)\lambda_{0},\dots,(k+n-3)\lambda_{0},l\lambda_{0}),
$$
$$
rang(\g)=1.
$$
\end{enumerate}

\item Si ${\rm gr}\frak{g}\simeq \frak{Q_{n-1}}\oplus \mathbb{C} \quad (n\ge 7, \: n \; {\rm impair})$ alors $ty(\g)=3$ et
\begin{enumerate}
\item $\frak{g}=Q_{n-1}\oplus \mathbb{C}$
$$
\begin{array}{lll}
\lbrack Y_{0},Y_{i}\rbrack=Y_{i+1}, & 1 \le i \le n-4,\\
\lbrack Y_{i},Y_{n-i-1}\rbrack=(-1)^{i-1}Y_{n-2}, & 1\leq i\leq \frac{n-3}{2},\\
\end{array}
$$
$$
f\sim diag(\lambda_{0},\lambda_{1},\lambda_{0}+\lambda_{1},2\lambda_{0}+\lambda_{1},\dots,(n-4)\lambda_{0}+\lambda_{1},(n-4)\lambda_{0}+2\lambda_{1},\lambda_{n-1}),
$$
$$
rang(\g)=3.
$$

\item \label{2.(b)} $\frak{g}=B_{n-1}^{k}(\alpha_{1},\dots,\alpha_{t-1})\oplus \mathbb{C} \quad t=[\frac{n-k-1}{2}] \quad 2 \le k \le n-5$
$$
\begin{array}{lll}
\lbrack Y_{0},Y_{i}\rbrack=Y_{i+1}, & 1 \le i \le n-4, &\\
\lbrack Y_{i},Y_{n-i-2}\rbrack =(-1)^{i-1}Y_{n-2}, & 1\leq i\leq \frac{n-3}{2}, &\\
\lbrack Y_{i},Y_{i+1}\rbrack=\alpha_{i} Y_{2i+k-1}, & 1 \le i \le t-1, &\\
\lbrack Y_{i},Y_{j}\rbrack=a_{i,j} Y_{i+j+k-1}, & 1 \le i <j, & i+j \le n-k-2,\\
\end{array}
$$
$$
f\sim diag(\lambda_{0},k\lambda_{0},(k+1)\lambda_{0},(k+2)\lambda_{0},\dots,(n-4+k)\lambda_{0},(n-4+2k)\lambda_{0},\lambda_{n-1}),
$$
$$
rang(\g)=2.
$$

\item \label{2.(c)} $\frak{g}=Q_{n-1}\overrightarrow{\oplus}_{l}^{a} \mathbb{C}\quad 2 \le l \le n-4$
$$
\begin{array}{ll}
\lbrack Y_{0},Y_{i}\rbrack=Y_{i+1}, & 1 \le i \le n-4,\\
\lbrack Y_{i},Y_{n-i-2}\rbrack=(-1)^{i-1}Y_{n-2}, & 1\leq i\leq \frac{n-3}{2},\\
\lbrack Y_{i},Y_{n-1}\rbrack=Y_{i+l}, & 1\leq i\leq n-l-3,\\
\end{array}
$$
$$
f\sim diag(\lambda_{0},\lambda_{1},\lambda_{0}+\lambda_{1},2\lambda_{0}+\lambda_{1},\dots,(n-4)\lambda_{0}+\lambda_{1},(n-4)\lambda_{0}+2\lambda_{1},l\lambda_{0}),
$$
$$
rang(\g)=2.
$$
\item \label{2.(d)} $\frak{g}=B_{n-1}^{k}(\alpha_{1},\dots,\alpha_{t-1})\overrightarrow{\oplus}_{l}^{a} \mathbb{C} \quad t=[\frac{n-k-1}{2}] \quad 2 \le k \le n-5, \; 2 \le l \le n-4$
$$
\begin{array}{lll}
\lbrack Y_{0},Y_{i}\rbrack=Y_{i+1}, & 1 \le i \le n-4, &\\
\lbrack Y_{i},Y_{n-i-2}\rbrack=(-1)^{i-1}Y_{n-2}, & 1\leq i\leq \frac{n-3}{2}, &\\
\lbrack Y_{i},Y_{i+1}\rbrack=\alpha_{i} Y_{2i+k}, & 1 \le i \le t-1, &\\
\lbrack Y_{i},Y_{j}\rbrack=a_{i,j} Y_{i+j+k-1}, & 1 \le i <j, & i+j \le n-k-2,\\
\lbrack Y_{i},Y_{n-1}\rbrack=Y_{i+l}, & 1\leq i\leq n-l-3, &\\
\end{array}
$$
$$
f\sim diag(\lambda_{0},k\lambda_{0},(k+1)\lambda_{0},(k+2)\lambda_{0},\dots,(n-4+k)\lambda_{0},(n-4+2k)\lambda_{0},l\lambda_{0}),
$$
$$
rang(\g)=1.
$$
\item \label{2.(e)} $\frak{g}=Q_{n-1}\overrightarrow{\oplus}_{l}^{b} \mathbb{C}\quad 2 \le l \le n-4$
$$
\begin{array}{ll}
\lbrack Y_{0},Y_{i}\rbrack=Y_{i+1}, & 1 \le i \le n-4,\\
\lbrack Y_{i},Y_{n-i-2}\rbrack=(-1)^{i-1}Y_{n-2}, & 1\leq i\leq \frac{n-3}{2},\\
\lbrack Y_{0},Y_{n-1}\rbrack=Y_{n-2}, &\\
\lbrack Y_{i},Y_{n-1}\rbrack=Y_{i+l}, & 1\leq i\leq n-l-3,\\
\end{array}
$$
$$
f\sim diag(\lambda_{0},\beta \lambda_{0},(k+1)\lambda_{0},(\beta+2)\lambda_{0},\dots,(n-4+\beta)\lambda_{0},(n-4+2\beta)\lambda_{0},(n-5+2\beta)\lambda_{0})
$$
o\`u $\beta=\frac{l-n+5}{2},$
$$
rang(\g)=1.
$$
\item \label{2.(f)} $\frak{g}=Q_{n-1}\overrightarrow{\oplus}^{c} \mathbb{C}$
$$
\begin{array}{ll}
\lbrack Y_{0},Y_{i}\rbrack=Y_{i+1}, & 1 \le i \le n-4,\\
\lbrack Y_{i},Y_{n-i-2}\rbrack=(-1)^{i-1}Y_{n-2}, & 1\leq i\leq \frac{n-3}{2},\\
\lbrack Y_{0},Y_{n-1}\rbrack=Y_{n-2}, &\\
\end{array}
$$
$$
f\sim diag(\lambda_{0},\lambda_{1},\lambda_{0}+\lambda_{1},2\lambda_{0}+\lambda_{1},\dots,(n-4)\lambda_{0}+\lambda_{1},(n-4)\lambda_{0}+2\lambda_{1},(n-5)\lambda_{0}+2\lambda_{1}),
$$
$$
rang(\g)=2.
$$
\item \label{2.(g)} $\frak{g}=B_{n-1}^{k}(\alpha_{1},\dots,\alpha_{t-1})\overrightarrow{\oplus}^{c} \mathbb{C} \quad t=[\frac{n-k-1}{2}] \quad 2 \le k \le n-5$
$$
\begin{array}{lll}
\lbrack Y_{0},Y_{i}\rbrack=Y_{i+1}, & 1 \le i \le n-4, &\\
\lbrack Y_{i},Y_{n-i-2}\rbrack=(-1)^{i-1}Y_{n-2}, & 1\leq i\leq \frac{n-3}{2}, &\\
\lbrack Y_{i},Y_{i+1}\rbrack=\alpha_{i} Y_{2i+k}, & 1 \le i \le t-1, &\\
\lbrack Y_{i},Y_{j}\rbrack=a_{i,j} Y_{i+j+k-1}, & 1 \le i <j, & i+j \le n-k-2,\\
\lbrack Y_{0},Y_{n-1}\rbrack=Y_{n-2}, & &\\
\end{array}
$$
$$
f\sim diag(\lambda_{0},k\lambda_{0},(k+1)\lambda_{0},(k+2)\lambda_{0},\dots,(n-4+k)\lambda_{0},(n-4+2k)\lambda_{0},(n-5+2k)\lambda_{0}),
$$
$$
rang(\g)=1.
$$
\end{enumerate}

\item Si ${\rm gr}\frak{g}\simeq \frak{L_{n,r}} \quad (n\ge5, \;r\,impair, \;3\le r\le 2[\frac{n-1}{2} ]-1)$ alors $ty(\g)=2$ et
\begin{enumerate}
\item \label{3.(a)} $\frak{g}=\frak{L_{n,r}}$
$$
\begin{array}{ll}
\lbrack Y_{0},Y_{i} \rbrack=Y_{i+1}, &  i=1,\dots,n-3,\\
\lbrack Y_{i},Y_{r-i} \rbrack=(-1)^{i-1}Y_{n-1}, &  i=1,\dots,\frac{r-1}{2}\\
\end{array}
$$
$$
f\sim diag(\lambda_{0},\lambda_{1},\lambda_{0}+\lambda_{1},2\lambda_{0}+\lambda_{1},\dots,(n-3)\lambda_{0}+\lambda_{1},(r-2)\lambda_{0}+2\lambda_{1}),
$$
$$
rang(\g)=2.
$$
\item \label{3.(b)} $\frak{g}=\frak{C_{n,r}^{k}}(\alpha_{1},\dots,\alpha_{t-1}), \quad 2\le k\le n-4,\;t=[\frac{n-k}{2}]$
$$
\begin{array}{ll}
\lbrack Y_{0},Y_{i} \rbrack=Y_{i+1}, &  i=1,\dots,n-3\\
\lbrack Y_{i},Y_{r-i} \rbrack=\left\{\begin{array}{l}(-1)^{i-1}Y_{n-1}+a_{i,r-i}Y_{r+k-1}\\ (-1)^{i-1}Y_{n-1}\\ \end{array}\right. &  \begin{array}{ll} si\: k\le n-r-1, & i=1,\dots,\frac{r-1}{2}\\ si\: k> n-r-1, & i=1,\dots,\frac{r-1}{2}\\ \end{array}\\
\lbrack Y_{i},Y_{i+1} \rbrack=\alpha_{i}Y_{2i+k}, &  i=1,\dots,t-1\\
\lbrack Y_{i},Y_{j} \rbrack=a_{i,j}Y_{i+j+k-1}, &  1\le i<j<n-1, \;r\ne i+j\le n-k-1,\\
\lbrack Y_{i},Y_{n-1} \rbrack=Y_{2k+r+i-2}, &  i=1,\dots,n-r-2k\\
\end{array}
$$
$$
f\sim diag(\lambda_{0},k\lambda_{0},(1+k)\lambda_{0},(2+k)\lambda_{0},\dots,(n-3+k)\lambda_{0},(r-2+2k)\lambda_{0}),
$$
$$
rang(\g)=1.
$$
\item \label{3.(c)} $\frak{g}=\frak{D_{n,r}^{k}}, \quad 1\le k\le [\frac{n-r-2}{2}]$
$$
\begin{array}{ll}
\lbrack Y_{0},Y_{i} \rbrack=Y_{i+1}, &  i=1,\dots,n-3\\
\lbrack Y_{i},Y_{r-i} \rbrack=(-1)^{i-1}Y_{n-1}, &  i=1,\dots,\frac{r-1}{2}\\
\lbrack Y_{i},Y_{n-1} \rbrack=Y_{2k+r+i-1}, &  i=1,\dots,n-r-2k-1\\
\end{array}
$$
$$
f\sim diag(\lambda_{0},(k+\frac{1}{2})\lambda_{0},(k+\frac{3}{2})\lambda_{0},(k+\frac{5}{2})\lambda_{0},\dots,(k+\frac{2n-5}{2})\lambda_{0},(r-1+2k)\lambda_{0})
$$
$$
rang(\g)=1.
$$
\end{enumerate}

\item Si ${\rm gr}\frak{g}\simeq \frak{Q_{n,r}}\quad (n\ge7,\; n\,impair, \;r\,impair, \;3\le r\le n-4)$ alors $ty(\g)=2$ et
\begin{enumerate}
\item \label{4.(a)} $\frak{g}=\frak{Q_{n,r}}$
$$
\begin{array}{ll}
\lbrack Y_{0},Y_{i} \rbrack=Y_{i+1}, &  i=1,\dots,n-4\\
\lbrack Y_{i},Y_{r-i} \rbrack=(-1)^{i-1}Y_{n-1}, &  i=1,\dots,\frac{r-1}{2}\\
\lbrack Y_{i},Y_{n-2-i} \rbrack=(-1)^{i-1}Y_{n-2}, &  i=1,\dots,\frac{n-3}{2}\\
\end{array}
$$
$$
f\sim diag(\lambda_{0},\lambda_{1},\lambda_{0}+\lambda_{1},2\lambda_{0}+\lambda_{1},\dots,(n-4)\lambda_{0}+\lambda_{1},(n-4)\lambda_{0}+2\lambda_{1},(r-2)\lambda_{0}+2\lambda_{1}),
$$
$$
rang(\g)=2.
$$
\item \label{4.(b)} $\frak{g}=\frak{E_{n,r}^{k}}(\alpha_{1},\dots,\alpha_{t-1}), \quad2\le k\le n-5,\;t=[\frac{n-k-1}{2}]$
$$
\begin{array}{ll}
\lbrack Y_{0},Y_{i} \rbrack=Y_{i+1}, &  i=1,\dots,n-4\\
\lbrack Y_{i},Y_{r-i} \rbrack=\left\{\begin{array}{l}(-1)^{i-1}Y_{n-1}+a_{i,r-i}Y_{r+k-1}\\ (-1)^{i-1}Y_{n-1}\\ \end{array}\right. &   \begin{array}{ll}si\, k\le n-r-2, & i=1,\dots,\frac{r-1}{2}\\ si\, k> n-r-2, & i=1,\dots,\frac{r-1}{2}\\ \end{array} \\
\lbrack Y_{i},Y_{n-2-i} \rbrack=(-1)^{i-1}Y_{n-2}, &  i=1,\dots,\frac{n-3}{2}\\
\lbrack Y_{i},Y_{i+1} \rbrack=\alpha_{i}Y_{2i+k}, &  i=1,\dots,t-1\\
\lbrack Y_{i},Y_{j} \rbrack=a_{i,j}Y_{i+j+k-1}, &  1\le i<j<n-1, \;r\ne i+j\le n-k-2\\
\lbrack Y_{i},Y_{n-1} \rbrack=Y_{2k+r+i-2}, &  i=1,\dots,n-r-2k-1\\
\end{array}
$$
$$
f\sim diag(\lambda_{0},k\lambda_{0},(1+k)\lambda_{0},(2+k)\lambda_{0},\dots,(n-4+k)\lambda_{0},(n-4+2k)\lambda_{0},(r-2+2k)\lambda_{0}),
$$
$$
rang(\g)=1.
$$
\item \label{4.(c)} $\frak{g}=\frak{F_{n,r}^{k}}, \quad 1\le k\le [\frac{n-r-4}{2}]$
$$
\begin{array}{ll}
\lbrack Y_{0},Y_{i} \rbrack=Y_{i+1}, &  i=1,\dots,n-4\\
\lbrack Y_{i},Y_{r-i} \rbrack=(-1)^{i-1}Y_{n-1}, &  i=1,\dots,\frac{r-1}{2}\\
\lbrack Y_{i},Y_{n-2-i} \rbrack=(-1)^{i-1}Y_{n-2}, &  i=1,\dots,\frac{n-3}{2}\\
\lbrack Y_{i},Y_{n-1} \rbrack=Y_{2k+r+i-1}, &  i=1,\dots,n-r-2k-2\\
\end{array}
$$
$$
f\sim diag(\lambda_{0},(k+\frac{1}{2})\lambda_{0},(k+\frac{3}{2})\lambda_{0},(k+\frac{5}{2})\lambda_{0},\dots,(k+\frac{2n-7}{2})\lambda_{0},(n+2k-3)\lambda_{0},(r+2k-1)\lambda_{0}),
$$
$$
rang(\g)=1.
$$
\end{enumerate}

\item Si ${\rm gr}\frak{g}\simeq \frak{T_{n,n-4}}\quad (n\ge7,\; n\,impair)$ alors $ty(\g)=2$ et
\begin{enumerate}
\item \label{5.(a)} $\frak{g}=\frak{T_{n,n-4}}$
$$
\begin{array}{ll}
\lbrack Y_{0},Y_{i} \rbrack=Y_{i+1}, &  i=1,\dots,n-5\\
\lbrack Y_{0},Y_{n-3} \rbrack=Y_{n-2}, &\\
\lbrack Y_{0},Y_{n-1} \rbrack=Y_{n-3}, &\\
\lbrack Y_{i},Y_{n-4-i} \rbrack=(-1)^{i-1}Y_{n-1}, &  i=1,\dots,\frac{n-5}{2}\\
\lbrack Y_{i},Y_{n-3-i} \rbrack=(-1)^{i-1} \frac{n-3-2i}{2} Y_{n-3}, &  i=1,\dots,\frac{n-5}{2}\\
\lbrack Y_{i},Y_{n-2-i} \rbrack=(-1)^{i} (i-1) \frac{n-3-i}{2} Y_{n-2},\quad &  i=2,\dots,\frac{n-3}{2}\\
\end{array}
$$
$$
f\sim diag(\lambda_{0},\lambda_{1},\lambda_{0}+\lambda_{1},2\lambda_{0}+\lambda_{1},\dots,(n-5)\lambda_{0}+\lambda_{1},(n-5)\lambda_{0}+2\lambda_{1},(n-4)\lambda_{0}+2\lambda_{1},(n-6)\lambda_{0}+2\lambda_{1}),
$$
$$
rang(\g)=2.
$$
\item \label{5.(b)} $\frak{g}=\frak{G_{n,r}^{k}}(\alpha_{1},\dots,\alpha_{t-1}),\quad 2\le k\le n-6,\;t=[\frac{n-k-2}{2}]$
$$
\begin{array}{ll}
\lbrack Y_{0},Y_{i} \rbrack=Y_{i+1}, &  i=1,\dots,n-5\\
\lbrack Y_{0},Y_{n-3} \rbrack=Y_{n-2},\\
\lbrack Y_{0},Y_{n-1} \rbrack=Y_{n-3},\\
\lbrack Y_{1},Y_{n-1} \rbrack=Y_{n-2} &si \; k=2\\
\lbrack Y_{i},Y_{n-4-i} \rbrack=(-1)^{i-1}Y_{n-1}, &  i=1,\dots,\frac{n-5}{2}\\
\lbrack Y_{i},Y_{n-3-i} \rbrack=(-1)^{i-1}\frac{n-3-2i}{2}Y_{n-3}, &  i=1,\dots,\frac{n-5}{2}\\
\lbrack Y_{i},Y_{n-2-i} \rbrack=(-1)^{i} (i-1) \frac{n-2-i}{2} Y_{n-2}, &  i=1,\dots,\frac{n-3}{2},\\
\lbrack Y_{i},Y_{i+1} \rbrack=\alpha_{i}Y_{2i+k}, &  i=1,\dots,t-1\\
\lbrack Y_{i},Y_{j} \rbrack=a_{i,j}Y_{i+j+k-1}, &  1\le i<j<n-2,\, i+j\le n-k-3\\
\end{array}
$$
$$
f\sim diag(\lambda_{0},k\lambda_{0},(1+k)\lambda_{0},(2+k)\lambda_{0},\dots,(n-5+k)\lambda_{0},(n-5+2k)\lambda_{0},(n-4+2k)\lambda_{0},(n-6+2k)\lambda_{0}),
$$
$$
rang(\g)=1.
$$
\end{enumerate}

\item Si ${\rm gr}\frak{g}\simeq \frak{T_{n,n-3}}\quad (n\ge6,\; n\,pair)$ alors $ty(\g)=2$ et
\begin{enumerate}
\item \label{6.(a)} $\frak{g}=\frak{T_{n,n-3}}$
$$
\begin{array}{ll}
\lbrack Y_{0},Y_{i} \rbrack=Y_{i+1}, &  i=1,\dots,n-4\\
\lbrack Y_{0},Y_{n-1} \rbrack=Y_{n-2},\\
\lbrack Y_{i},Y_{n-3-i} \rbrack=(-1)^{i-1}Y_{n-1}, &  i=1,\dots,\frac{n-4}{2}\\
\lbrack Y_{i},Y_{n-2-i} \rbrack=(-1)^{i-1} \frac{n-2-2i}{2} Y_{n-2}, &  i=1,\dots,\frac{n-4}{2}\\
\end{array}
$$
$$
f\sim diag(\lambda_{0},\lambda_{1},\lambda_{0}+\lambda_{1},2\lambda_{0}+\lambda_{1},\dots,(n-4)\lambda_{0}+\lambda_{1},(n-4)\lambda_{0}+2\lambda_{1},(n-5)\lambda_{0}+2\lambda_{1}),
$$
$$
rang(\g)=2.
$$
\item \label{6.(b)} $\frak{g}=\frak{H_{n,r}^{k}}(\alpha_{1},\dots,\alpha_{t-1}), \quad2\le k\le n-5,\;t=[\frac{n-k-1}{2}]$
$$
\begin{array}{ll}
\lbrack Y_{0},Y_{i} \rbrack=Y_{i+1}, &  i=1,\dots,n-4\\
\lbrack Y_{0},Y_{n-1} \rbrack=Y_{n-2},\\
\lbrack Y_{i},Y_{n-3-i} \rbrack=(-1)^{i-1}Y_{n-1}, &  i=1,\dots,\frac{n-4}{2}\\
\lbrack Y_{i},Y_{n-2-i} \rbrack=(-1)^{i-1} \frac{n-2-2i}{2} Y_{n-2}, &  i=1,\dots,\frac{n-4}{2}\\
\lbrack Y_{i},Y_{i+1} \rbrack=\alpha_{i}Y_{2i+k}, &  i=1,\dots,t-1\\
\lbrack Y_{i},Y_{j} \rbrack=a_{i,j}Y_{i+j+k-1}, &  1\le i<j<n-2,\, i+j\le n-k-2\\
\end{array}
$$
$$
f\sim diag(\lambda_{0},k\lambda_{0},(1+k)\lambda_{0},(2+k)\lambda_{0},\dots,(n-4+k)\lambda_{0},(n-4+2k)\lambda_{0},(n-5+2k)\lambda_{0}),
$$
$$
rang(\g)=1.
$$
\end{enumerate}

\item Si ${\rm gr}\frak{g}\simeq \frak{E_{9,5}^{1}}$ alors $\frak{g}\simeq \frak{E_{9,5}^{1}}$ et $ty(\g)=2$
$$
\begin{array}{llllll}
& \lbrack Y_{0},Y_{i}\rbrack =Y_{i+1}, & i=1,\dots,4, & \lbrack
Y_{0},Y_{8}\rbrack =Y_{6}, & \lbrack Y_{1},Y_{4}\rbrack =Y_{8}, &
\lbrack Y_{1},Y_{5}\rbrack =2Y_{6},\\
& \lbrack Y_{1},Y_{6}\rbrack =3Y_{7}, & \lbrack Y_{2},Y_{3}\rbrack =-Y_{8}, &
 \lbrack Y_{2},Y_{4}\rbrack =-Y_{6}, & \lbrack Y_{2},Y_{8}\rbrack =-3Y_{7},\\
\end{array}
$$
$$
f\sim
diag(\lambda_{0},\lambda_{1},\lambda_{0}+\lambda_{1},2\lambda_{0}+\lambda_{1},3\lambda_{0}+\lambda_{1},4\lambda_{0}+\lambda_{1},4\lambda_{0}+2\lambda_{1},4\lambda_{0}+3\lambda_{1},3\lambda_{0}+2\lambda_{1}),
$$
$$
rang(\g)=2.
$$

\item Si ${\rm gr}\frak{g}\simeq \frak{E_{9,5}^{2}}$ alors $\frak{g}\simeq \frak{E_{9,5}^{2}}$ et $ty(\g)=2$
$$
\begin{array}{lll}
\lbrack Y_{0},Y_{i}\rbrack =Y_{i+1}, &i=1,\dots,6 & \lbrack Y_{0},Y_{8}\rbrack =Y_{6},\\
\lbrack Y_{1},Y_{4}\rbrack =Y_{8}, & \lbrack Y_{1},Y_{5}\rbrack =2Y_{6}, & \lbrack Y_{1},Y_{6}\rbrack =Y_{7},\\
\lbrack Y_{2},Y_{3}\rbrack =-Y_{8}, & \lbrack Y_{2},Y_{4}\rbrack =-Y_{6}, & \lbrack Y_{2},Y_{5}\rbrack =Y_{7},\\
\lbrack Y_{2},Y_{8}\rbrack =-Y_{7}, & \lbrack Y_{3},Y_{4}\rbrack =-2Y_{7}, &\\
\end{array}
$$
$$
f\sim diag(\lambda_{0},\lambda_{0},2\lambda_{0},3\lambda_{0},4\lambda_{0},5\lambda_{0},6\lambda_{0},7\lambda_{0},5\lambda_{0}),
$$
$$
rang(\g)=1.
$$

\item Si ${\rm gr}\frak{g}\simeq \frak{E_{9,5}^{3}}$ alors $\frak{g}\simeq \frak{E_{9,5}^{3}}$ et $ty(\g)=2$
$$
\begin{array}{llllll}
 & \lbrack Y_{0},Y_{i}\rbrack =Y_{i+1}, & i=1,\dots,4, & \lbrack Y_{0},Y_{6}\rbrack =Y_{7}, & \lbrack Y_{0},Y_{8}\rbrack =Y_{6}, & \lbrack Y_{1},Y_{4}\rbrack =Y_{8},\\
& \lbrack Y_{2},Y_{3}\rbrack =-Y_{8}, & \lbrack Y_{2},Y_{4}\rbrack =-Y_{6}, & \lbrack Y_{2},Y_{5}\rbrack =2Y_{7}, & \lbrack Y_{3},Y_{4}\rbrack =-3Y_{7},\\
\end{array}
$$
$$
f\sim
diag(\lambda_{0},\lambda_{1},\lambda_{0}+\lambda_{1},2\lambda_{0}+\lambda_{1},3\lambda_{0}+\lambda_{1},4\lambda_{0}+\lambda_{1},4\lambda_{0}+2\lambda_{1},5\lambda_{0}+2\lambda_{1},3\lambda_{0}+2\lambda_{1}).
$$
$$
rang(\g)=2.
$$

\item Si ${\rm gr}\frak{g}\simeq \frak{E_{7,3}}$ alors $\frak{g}\simeq \frak{E_{7,3}}$ et $ty(\g)=2$
$$
\begin{array}{lll}
\lbrack Y_{0},Y_{i}\rbrack =Y_{i+1}, & i=1,\dots,4 & \lbrack Y_{0},Y_{6}\rbrack =Y_{4},\\
\lbrack Y_{1},Y_{2}\rbrack =Y_{6}, & \lbrack Y_{1},Y_{3}\rbrack =Y_{4}, & \lbrack Y_{1},Y_{4}\rbrack =Y_{5},\\
\lbrack Y_{2},Y_{6}\rbrack =-Y_{5},& &
\end{array}\\
$$
$$
f\sim diag(\lambda_{0},\lambda_{0},2\lambda_{0},3\lambda_{0},4\lambda_{0},5\lambda_{0},3\lambda_{0}),
$$
$$
rang(\g)=1.
$$
\end{enumerate}
Les param\`{e}tres $(\alpha_{1},\dots,\alpha_{t-1})$ v\'{e}rifient les relations polynomiales d\'{e}coulant des identit\'{e}s de Jacobi et les constantes $a_{i,j}$ v\'erifient le syst\`{e}me:
$$
\begin{array}{l}
a_{i,i}=0,\\
a_{i,i+1}=\alpha_{i},\\
a_{i,j}=a_{i+1,j}+a_{i,j+1}.
\end{array}
$$
\end{theo}

\begin{coro}
Soit $\frak{g}$ une \al quasi-filiforme de rang non nul. L'alg\`ebre $\frak{g}$ est de rang maximal si et seulement si elle est isomorphe \`a une des alg\`ebres suivantes: $\frak{g}=L_{n-1} \overrightarrow{\oplus} \mathbb{C}$, $\frak{g}=Q_{n-1} \overrightarrow{\oplus} \mathbb{C}$, $\frak{L_{n,r}}$, $\frak{Q_{n,r}}$, $\frak{T_{n,n-4}}$, $\frak{T_{n,n-3}}$, $\frak{E_{9,5}^{1}}$ ou bien $\frak{E_{9,5}^{3}}$.
\end{coro}

Les deux th\'eor\`emes ci-dessous nous donnent des conditions
suffisantes pour la compl\'etude d'une alg\`ebre de Lie.

\begin{theo}
\label{teo1}
Soit $\frak{n}$ une \al nilpotente de rang maximal et $\frak{h}$ un tore maximal de $\frak{n}$. Alors l'alg\`ebre $\frak{g}=\frak{h}\overrightarrow{\oplus}\frak{n}$ est compl\`ete.
\end{theo}
\begin{theo}
\label{teo2}
Soient $\frak{g}$ une \al et $\frak{h}$ une sous-alg\`ebre de Cartan qui v\'erifient les conditions suivantes:
\begin{enumerate}
\item  $\frak{h}$ est ab\'elienne.

\item  $\frak{g}$ se d\'ecompose de la fa\c{c}on $\frak{h}\oplus \sum_{\alpha \in \Delta }\frak{g}_{\alpha }$ avec $\Delta \subset \frak{h}^{\ast }-\{0\}$.

\item  Il existe un syst\`eme de g\'en\'erateurs $\{ \alpha _{1},..,\alpha_{l}\}\subseteq\Delta$ de $\frak{h}^{\ast}$ tel que $\dim \frak{g}_{\alpha _{j}}=1$ pour $1\le j\le l$.

\item Soit $\{\alpha _{1},..,\alpha _{r}\} $ une base de $\frak{h}^{\ast }$. Pour $r+1\leq s\leq l,$ on a:
\begin{equation*}
\alpha _{s}=\sum_{i=1}^{t}k_{is}\alpha _{j_{i}}-\sum_{i=1+t}^{r}k_{is}\alpha_{j_{i}}
\end{equation*}
o\`u  $k_{is}\in \mathbb{N\cup \{}0\},\left( j_{1},..,j_{r}\right) $ est une permutation de $\left( 1,..,r\right) $,et il existe une formule
\begin{eqnarray*}
&&[\underset{k_{1s}}{\underbrace{x_{j_{1}},..,x_{j_{1}}}},..\underset{k_{ts}%
}{\underbrace{x_{j_{t}},..,x_{j_{t}}},}...,x_{k_{m}}] \\
&=&[\underset{k_{t+1s}}{\underbrace{x_{j_{t}+1},..,x_{j_{t}+1}}},..\underset{%
k_{rs}}{\underbrace{x_{j_{r}},..,x_{j_{r}}},x_{s},x_{k_{1}}}...,x_{k_{m}}]
\end{eqnarray*}
 l'ordre de calcul des crochets n'ayant pas d'importance, $0\neq x_{j}\in \frak{g}_{\alpha _{j}}$ et $m\neq 0$ si $t=r$.
\end{enumerate}
Alors $\frak{g}$ est une \al compl\`ete.
\end{theo}
Ces th\'eor\`emes sont d\'emontr\'es dans les articles ~\cite{Zhu} et ~\cite{Meng}.
\begin{theo}
Toute \al quasi-filiforme de rang non nul est compl\'etable.
\end{theo}
\dd
Soit $\frak{n}$ une \al quasi-filiforme de rang non nul, elle est donc isomorphe \`a une des alg\`ebres du th\'eor\`eme \ref{theo}. On consid\`ere la somme semi-directe $\frak{g}=\frak{h}\overrightarrow{\oplus}\frak{n}$ o\`u $\frak{h}$ est un tore maximal associ\'e \`a $\frak{n}$.\\
Quand $\frak{n}$ est isomorphe \`a $\frak{E_{9,5}^{2}}$ ou bien \`a $\frak{E_{7,3}}$, par le calcul des deux premiers groupes de cohomologie, on v\'erifie que $\g$ est compl\`ete et donc $\n$ est compl\'etable.\\
Si $\frak{n}$ est isomorphe \`a $L_{n-1} \oplus \mathbb{C}$, $Q_{n-1} \oplus \mathbb{C}$, $\frak{L_{n,r}}$, $\frak{Q_{n,r}}$, $\frak{T_{n,n-4}}$, $\frak{T_{n,n-3}}$, $\frak{E_{9,5}^{1}}$ ou bien \`a $\frak{E_{9,5}^{3}}$, $\frak{n}$ est de rang maximal et d'apr\`es le th\'eor\`eme \ref{teo1}, $\frak{g}$ est compl\`ete.\\
Si $\frak{n}\simeq A_{n-1}^{k}(\alpha_{1},\dots,\alpha_{t-1})\oplus \mathbb{C}$, $\g$ se d\'ecompose alors en somme directe $\g=\g_{1}\oplus\g_{2}$ o\`u $\g_{1}=A_{n-1}^{k}(\alpha_{1},\dots,\alpha_{t-1})\overrightarrow{\oplus} \w$ , $\w$ \'etant un tore maximal de $A_{n-1}^{k}(\alpha_{1},\dots,\alpha_{t-1})$ et $\g_{2}$ l'alg\`ebre non-ab\'elienne de dimension $2$. Dans ~\cite{AnCa}, on d\'emontre que $\g_{1}$ est compl\`ete et comme $H^{0}(\g_{2},\g_{2})=H^{1}(\g_{2},\g_{2})=0$, il en r\'esulte que $\g$ est compl\`ete.\\
De fa\c{c}on analogue, on prouve la compl\'etude de $\g$ lorsque $\n$ est isomorphe \`a $B_{n-1}^{k}(\alpha_{1},\dots,\alpha_{t-1})\oplus \mathbb{C}$.\\
Si $\frak{n}\simeq \frak{D_{n,r}^{k}}$, $\Delta=\{ \alpha_{0}=2\lambda _{0},(\alpha_{j}=(2k+2j-1)\lambda _{0})_{1\leq j \leq n-2},\alpha_{n-1}=2(r-1+2k)\lambda _{0}\}$ est un syst\`eme de poids de $\n$. Les deux premi\`eres conditions du th\'eor\`eme \ref{teo2} se v\'erifient et d'apr\`es le lemme 2.3 de ~\cite{Meng}
$$
Der (\g) = D_{0}+ ad (\g)
$$
o\`u $D_{0}=\{ \phi \in Der (\g) / \phi(h)=0 \forall h\in \frak{h}\}$. Pour prouver que $\g$ est compl\`ete, il suffit de voir que $D_{0}\subseteq ad (\g)$. Soit $D\in D_{0}$, pour tout $X_{i}\in\g_{\alpha_{i}}(0\leq i\leq n-1)$ et $h\in\frak{h}$, on a:
$$
[h,D(X_{i})]=[h,D(X_{i})]+[D(h),X_{i}]=D([h,X_{i}])=\alpha_{i}(h)D(X_{i}).
$$
Comme dim$\g_{\alpha_{i}}=1$, $D(Y_{i})=d_{i}Y_{i}$ pour $0\leq i\leq n-1$. A partir des crochets
$$
\begin{array}{ll}
\lbrack Y_{0},Y_{i} \rbrack=Y_{i+1}, &  i=1,\dots,n-3\\
\lbrack Y_{1},Y_{n-1} \rbrack=Y_{2k+r}, & \lbrack Y_{1},Y_{r-1} \rbrack=Y_{n-1},\\
\end{array}
$$
on obtient les relations
$$
d_{i}=d_{0}+(i-1)d_{1} \,\forall i\in\{1,\dots,n-2\},\,d_{n-1}=(2k+r-1)d_{0},\,2d_{1}=(2k+1)d_{0}.
$$
On en d\'eduit que $D$ est une d\'erivation int\'erieure de $\g$. De m\^eme, on d\'emontre que $\g$ est compl\`ete lorsque $\frak{n}\simeq \frak{F_{n,r}^{k}}$.\\
Supposons maintenant que $\frak{n}\simeq L_{n-1} \overrightarrow{\oplus}_{l} \mathbb{C}$, $\Delta=\{ \lambda _{0},\lambda _{1},(j\lambda _{0}+\lambda _{1})_{0\leq j \leq n-3},l\lambda _{0}\}$ est un syst\`eme de poids de $\n$. Les trois premi\`eres conditions du th\'eor\`eme \ref{teo2} se v\'erifient clairement, quant \`a la quatri\`eme, il suffit de remarquer que:
$$
[\overset{j}{\overbrace{Y_{0},[Y_{0},\dots[Y_{0},}}Y_{1}]\dots]]=Y_{j+1} \quad 1\leq j\leq n-3\\
$$
$$
[\overset{l}{\overbrace{Y_{0},[Y_{0},\dots[Y_{0},}}Y_{i}]\dots]]=Y_{i+l}=[-Y_{n-1},Y_{i}] \quad 1\leq i\leq n-l
$$
Ainsi, le th\'eor\`eme \ref{teo2} nous donne la compl\'etude de $\g$.\\
Pour les cas restants, on d\'emontre de la m\^eme fa\c{c}on que $\g$ est compl\`ete et on en conclut que $\frak{n}$ est compl\'etable.\\

\section{Sur le deuxi\`eme groupe de cohomologie des \als ayant un nilradical quasi-filiforme}
Nous avons vu que pour toute \al quasi-filiforme $\n$, l'alg\`ebre $\frak{g}=\frak{h}\overrightarrow{\oplus}\frak{n}$, $\frak{h}$ \'etant un tore maximal associ\'e \`a $\frak{n}$, est compl\`ete et donc $H^{0}(\g,\g)=H^{1}(\g,\g)=\{0\}$. Nous nous demandons ce qu'il en est du deuxi\`eme groupe de coholomogie.\\
Comme dans ~\cite{AnCa}, nous consid\'erons la famille d'alg\`ebres $\n =A_{n-1}^{k}(\alpha_{1},\dots,\alpha_{t-1})\oplus \mathbb{C}$ avec $t=[\frac{n-k}{2}]$ et $2 \le k \le n-4$ d\'efinie par
$$
\begin{array}{lll}
\lbrack Y_{0},Y_{i}\rbrack=Y_{i+1}, & 1 \le i \le n-3, &\\
\lbrack Y_{i},Y_{i+1}\rbrack=\alpha_{i} Y_{2i+k}, & 1 \le i \le t-1, &\\
\lbrack Y_{i},Y_{j}\rbrack=a_{i,j} Y_{i+j+k-1}, & 1 \le i <j, & i+j \le n-k-1\\
\end{array}
$$
o\`u
\begin{equation}
\label{Condaij}
\begin{array}{lll}
a_{i,i}=0, & a_{i,i+1}=\alpha_{i}, & a_{i,j}=a_{i+1,j}+a_{i,j+1}.
\end{array}
\end{equation}
Si $\lambda_{1}$ est non-nul, on peut prendre $\lambda_{1}=1$ et les param\`etres $\lambda_{2},\dots,\lambda_{t-1}$ distinguent les classes d'isomorphisme de la famille. Pour chaque $l\in\{2,\dots,t-1\}$, l'alg\`ebre $A_{n-1}^{k}(1,0,\dots,\alpha_{l},\dots,0)$ est une d\'eformation lin\'eaire de $A_{n-1}^{k}(1,0,\dots,0)$ qui correspond au $2$-cocycle suivant:
$$
F^{l}_{1}(X_{i},X_{j})=\gamma^{l}_{i,j}X_{i+j+k-1}
$$
Les coefficients $\gamma^{l}_{i,j}$ sont donn\'es par les relations $a_{i,j}=\sum_{l=1}^{t-1}\gamma^{l}_{i,j}\alpha_{l}$ qui d\'ecoulent des \'equations (\ref{Condaij}).\\
Toutes ces d\'eformations \'etant non-\'equivalentes, on trouve $t-2$ \'el\'ements de $Z^{2}(\g,\g)$ non-\'equivalents modulo $B^{2}(\g,\g)$ o\`u $\frak{g}=\frak{h}\overrightarrow{\oplus}(A_{n-1}^{k}(1,0,\dots,0)\oplus \mathbb{C})$. On en conclut que $\textrm{dim}H^{2}(\g,\g)\geq t-2$.
\begin{propo}
Pour tout $m\in \mathbb{N}^{+}$, il existe une \al $\frak{g}$ compl\`ete dont le nilradical est quasi-filiforme et telle que $\textrm{dim}H^{2}\left( \frak{g},\frak{g}\right) \geq m$.
\end{propo}

\noindent{\bf Remerciements:}
L'auteur est soutenu par les projets de recherche MTM2006-09152 et CCG07-UCM/ESP-2922, et remercie aussi la Fundaci\'{o}n Ram\'{o}n Areces qui finance sa bourse pr\'{e}doctorale.


\begin{thebibliography}{99}
\bibitem{AnCa}J. M. Ancochea Berm\'udez, R. Campoamor, Completable filiform Lie algebras, Linear Algebra and Appl. 367, 185-191, 2003.
\bibitem{Ca} R. Carles, Sur les suites d'alg\`ebres de Lie de d\'erivations, Arch. Math. 70, 262-269, 1998.
\bibitem{Fa} G. Favre, Syst\`eme de poids sur une alg\`ebre de Lie nilpotente, Manuscripta Math. 9, 53-90, 1973.
\bibitem{MFa} M. Favre, Alg\`ebres de Lie compl\`etes, C.R. Acad. Sci. Patris S\`er. I Math. 274, 1533-1535, 1972.
\bibitem{Ga} L. Garc\'{\i}a Vergnolle, Sur les alg\`ebres de Lie quasi-filiformes admettant un tore de d\'erivations. Manuscripta Math. 124, 489-505, 2007.
\bibitem{Gom} J.R. G\'omez, A. Jim\'enez-Merch\'an, Naturally graded quasi-filiform Lie algebras, J. Algebra 256, 221-228, 2002.
\bibitem{Go} M. Goze, Y. Hakimjanov, Sur les alg\`{e}bres de Lie nilpotentes admettant un tore de d\'{e}rivations, Manuscripta Math. 84, 115-124, 1994.
\bibitem{Ver} M. Vergne, Cohomologies des alg\`ebres de Lie nilpotentes. Application  \`a l'\'etude  de la vari\'et\'e des alg\`ebres de Lie nilpotentes, Bull. Soc. Math. France 98, 81-116, 1970.
\bibitem{Zhu} Zhu Lin Sheng, Daoji Meng, Solvable complete Lie algebras I. Communications in algebra, 24(13) 4181-4197, 1996.
\bibitem{Meng} Zhu Lin Sheng, Daoji Meng, Solvable complete Lie algebras II. Algebra Colloquium 289-296, 1998.
\bibitem{Vi} A. L. Onishchik, E. B. Vinberg, Lie Groups and Lie Algebras II, Springer Verlag, 1994.
\bibitem{Sch}  E. Schenkman, A theory of subinvariant Lie algebras, Amer. J. Math., 73, 453-474, 1951.



\end{thebibliography}
\end{document}